\documentclass[12pt]{amsart}
\usepackage{amssymb,latexsym}
\usepackage{booktabs,enumerate}

\makeatletter
\@namedef{subjclassname@2010}{%
  \textup{2010} Mathematics Subject Classification}
\makeatother

\newtheorem{theorem}{Theorem}[section]
\newtheorem{lemma}[theorem]{Lemma}

\numberwithin{equation}{section}

\DeclareMathOperator{\re}{Re}

\begin{document}

\baselineskip=17pt

\title{On sums of Ramanujan sums}

\author{T.H. Chan}
\address{Department of Mathematical Sciences, University of Memphis, Memphis, TN 38152, U.S.A.}
\email{tchan@memphis.edu}

\author{A.V. Kumchev}
\address{Department of Mathematics, 7800 York Road, Towson University, Towson, MD 21252, U.S.A.}
\email{akumchev@towson.edu}

\date{}

\begin{abstract}
  Let $c_q(n)$ denote the Ramanujan sum modulo $q$, and let $x$ and $y$ be large reals, with $x = o(y)$. We obtain asymptotic formulas for the sums 
  \[
    \sum_{n \le y} \Big( \sum_{q \le x} c_q(n) \Big)^k \quad (k = 1, 2).
  \]
\end{abstract}



\maketitle

\section{Introduction}

In this note, we study the moments of the average of the Ramanujan sum. To be precise, we want to evaluate the sums
\[
  C_k(x,y) = \sum_{n \le y} \Big( \sum_{q \le x} c_q(n) \Big)^k,
\]
where $k$ is a positive integer, $x$ and $y$ are large real numbers, and
\begin{equation}\label{eq0}
  c_q(n) = \sum_{\substack{1 \le a \le q\\ (a, q) = 1}} e(-an/q) = \sum_{d \mid  (q,n)} d\mu(q/d)
\end{equation}
is the Ramanujan sum. Our interest in this question stems from an old attempt \cite{ChKu05} to apply Fourier techniques to Diophantine approximations of reals by sums of rational numbers. At the time, we were surprised that information on the asymptotic behavior of $C_k(x,y)$ appeared to be missing from the literature. Of course, it is easy to evaluate $C_k(x,y)$ in some cases. For example, elementary arguments based on the second formula in \eqref{eq0} yield the asymptotic formulas
\[
   C_1(x,y) = y + O(x^{2})
\]
and
\begin{equation}\label{eq1}
   C_2(x,y) = \frac {yx^2}{2\zeta(2)} + O(x^4 + xy\log x),
\end{equation}
for any fixed $\epsilon > 0$. However, these bounds (as well as other simple things we have tried) are only of interest when $yx^{-2} \to \infty$ as $x \to \infty$. In our application to Diophantine approximations, on the other hand, we were interested in $C_2(x,y)$ when $x^{1+\epsilon} \le y \le x^{2 + \epsilon}$ for some fixed $\epsilon > 0$. Some numerical experimentation suggested that in that range the order of $C_2(x,y)$ is still $yx^2$, though the coefficient $\frac 1{2\zeta(2)}$ seemed somewhat off when $y$ is close to $x^2$. In the present note, we prove, among other things, that those empirical observations are true. Our first result is an asymptotic formula for $C_1(x,y)$.

\begin{theorem}\label{th1}
  Let $x$ be a large real number and $y \ge x$. Then
  \[
    C_1(x,y) = y - \frac{x^2}{4\zeta(2)} + O \big( xy^{1/3}\log x + x^3y^{-1} \big).
  \]
\end{theorem}

The proof of Theorem \ref{th1} uses exponential sum estimates and can be easily modified to yield different estimates for the remainder. The remainder claimed in the statement of the theorem is the result of applying one of the simplest exponential sum bounds---the exponential pair $(1/2, 1/2)$ from van der Corput's method \cite{GrKo91}. More sophisticated exponential pairs will result in replacing the error term $xy^{1/3}$ by different ones, which will be sharper for certain choices of the relative sizes of $x$ and $y$. The interested reader will have no trouble obtaining such improvements by appropriate modifications of \eqref{e6} below. Since our main focus is on maximizing the ``support'' of the result, we prefer the simple version above. 

Theorem \ref{th1} demonstrates that when $x = o(y)$, the Ramanujan sum $c_q(n)$ is $o(1)$ on average over $q \le x$ and $n \le y$. Our next theorem shows that this is no longer the case when the average $C_1(x,y)$ is replaced by the mean square sum $C_2(x,y)$.

\begin{theorem}\label{th2}
  Let $x$ be a large real number, $y \ge x$, and $B > 0$ be fixed.
  \begin{enumerate}[\upshape (i)]
    \item If $y \ge x^2(\log x)^B$, then \eqref{eq1} holds.
    \item If $x \le y \le x^2(\log x)^B$, then
    \[
      C_2(x,y) = \frac {yx^2}{2\zeta(2)}( 1 + 2\kappa(u) ) + O\big( yx^2(\log x)^{10}\big( x^{-1/2} + (y/x)^{-1/2} \big) \big),
    \]
    where $u = \log(yx^{-2})$, and $\kappa(u)$ is defined by \eqref{eq22} below and satisfies the inequalities
    \[
      \kappa(u) > -0.4, \qquad \kappa(u) \ll \exp( -|u|^{3/5 - \epsilon} )
    \]
    for any fixed $\epsilon >0$. In particular, $\kappa(u) = o(1)$ as $|u| \to \infty$.
  \end{enumerate}
\end{theorem}

The lower bound for $\kappa(u)$ is far from sharp and has been chosen so as to simplify the arguments in \S\ref{s5}. In fact, by increasing the amount of numeric calculations in \S\ref{s5} by an order of magnitude or two, it can be shown that $\kappa(u) > -0.3$. We chose not to pursue such a sharper bound here for the sake of clarity. However, for the benefit of the reader, we should mention that the more elaborate numerical computations that we carried suggest that $\kappa(u)$ is decreasing when $u < u_0$ and increasing when $u > u_0$, where $u_0 \approx 1.63$; its minimum value is approximately $\kappa(1.63) \approx -0.2943$. It also seems that $\kappa(u)$ is always negative.

It is natural to ask whether similar asymptotic formulas hold for $C_k(x,y)$ when $k \ge 3$. It appears that the analytic method in \S\ref{s4} below should yield some result, at least when $k=3$. However, since the technical details will most likely be quite ungainly, we leave such considerations to future work.

\section{The first moment}
\label{s2}

By \eqref{eq1}, we have
\[
  C_1(x,y) = \sum_{n \le y} \sum_{q \le x} \sum_{\substack{ d \mid q\\ d \mid n}} d \mu(q/d) 
  = \sum_{n \le y} \sum_{\substack{ dk \le x\\ d \mid n}} d \mu(k).
\]
Interchanging order of summation, we obtain
\begin{align}\label{e2}
  C_1(x,y) &= \sum_{dk \le x} d \mu(k) \sum_{\substack{ n \le y\\ d \mid n}} 1 
  = \sum_{dk \le x} d\mu(k) \left[ \frac yd \right] \\
  &= y \sum_{dk \le x} \mu(k) - \frac 12 \sum_{dk \le x} d \mu(k) - \sum_{dk \le x} d \mu(k)\psi(y/d) \notag\\
  &= C_{1,1}(x,y) - C_{1,2}(x,y) - C_{1,3}(x,y), \qquad \text{say.} \notag
\end{align}
Here, for a real $t$, $[t]$ is the integral part of $t$ and $\psi(t) = t - [t] - \frac 12$ is the saw-tooth function. It is easy to see that
\begin{equation}\label{e3}
  C_{1,1}(x,y) = y \sum_{dk \le x} \mu(k) = y \sum_{m \le x} \sum_{k \mid m} \mu(k) = y.
\end{equation}
Further,
\begin{align}\label{e4}
  C_{1,2}(x,y) &= \frac 12\sum_{dk \le x} d \mu(k) = \frac 12 \sum_{k \le x} \mu(k) \sum_{d \le x/k} d \\
  &= \frac 12 \sum_{k \le x} \mu(k) \left( \frac {x^2}{2k^2} + O\bigg( \frac xk \bigg) \right)  \notag\\
  &= \frac {x^2}4 \sum_{k \le x} \frac {\mu(k)}{k^2} + O(x\log x) = \frac {x^2}{4\zeta(2)} + O(x\log x). \notag
\end{align}

It remains to estimate $C_{1,3}(x,y)$. By partial summation,
\begin{equation}\label{e5}
  C_{1,3}(x,y) \ll \sum_{k \le x} \sum_{j = 0}^{\infty} N_{j} \sup_{\mathbf I} \Big| \sum_{d \in \mathbf I} \psi(y/d) \Big|,
\end{equation}
where $N_j = N_{j,k} = (x/k)2^{-j}$ and the supremum is over all subintervals $\mathbf I$ of $(N_{j}, 2N_{j}]$. We remark that the sum over $j$ is, in fact, finite and has $O(\log x)$ terms. We now use the following lemma, which is a special case of \cite[Lemma~4.3]{GrKo91} (see \cite{GrKo91} or \cite{Titc86} for the definition of exponential pairs).

\begin{lemma}\label{l3}
  Suppose that $(\varkappa,\lambda)$ is an exponential pair and that $\mathbf I$ is a subinterval of $(N, 2N]$. Then
  \[
    \sum_{n \in \mathbf I} \psi(y/n) \ll y^{\varkappa/(\varkappa+1)}N^{(\lambda-\varkappa)/(\varkappa+1)} + N^2y^{-1}.
  \]
\end{lemma}

We apply Lemma \ref{l3} with the exponential pair $(\varkappa, \lambda) = (1/2, 1/2)$ to the sum over $d$ on the right side of \eqref{e5}. This yields
\begin{align}\label{e6}
  C_{1,3}(x,y) &\ll \sum_{k \le x} \sum_{j = 0}^{\infty} ( N_{j}y^{1/3} + N_{j}^3y^{-1} ) 
  \ll \sum_{k \le x} ( (x/k)y^{1/3} + (x/k)^3y^{-1} ) \\
  &\ll ( xy^{1/3}\log x + x^3y^{-1} ). \notag
\end{align}
Theorem \ref{th1} now follows easily from \eqref{e2}--\eqref{e6}.

\section{The second moment: An elementary approach}
\label{s3}

In this section, we prove \eqref{eq1}. Similarly to \S\ref{s2}, we obtain
\[
  C_2(x,y) = \sum_{n \le y} \Big( \sum_{\substack{ dk \le x\\ d \mid n}} d\mu(k) \Big)^2
  = \sum_{d_1k_1 \le x} \sum_{d_2k_2 \le x} d_1d_2\mu(k_1)\mu(k_2) \left[ \frac {y}{[d_1,d_2]} \right],
\]
where $[d_1,d_2]$ denotes the least common multiple of $d_1$ and $d_2$. Hence,
\begin{equation}\label{eq3.1}
  C_2(x,y) = y \sum_{d_1k_1 \le x} \sum_{d_2k_2 \le x} (d_1,d_2)\mu(k_1)\mu(k_2) + O(E),
\end{equation}
where
\begin{equation}\label{eq3.2}
  E = y \sum_{d_1k_1 \le x} \sum_{d_2k_2 \le x} d_1d_2 \ll x^4.
\end{equation}
Furthermore, 
\begin{align}\label{eq3.3}
  \sum_{d_1k_1 \le x} \sum_{d_2k_2 \le x} (d_1,d_2)\mu(k_1)\mu(k_2) 
  &= \sum_{d \le x} d \mathop{\sum_{dl_1k_1 \le x} \sum_{dl_2k_2 \le x}}_{(l_1,l_2)=1} \mu(k_1)\mu(k_2) \\
  &= \sum_{d \le x} d \sum_{dl_1k_1 \le x} \sum_{dl_2k_2 \le x} \mu(k_1)\mu(k_2) 
  \sum_{l\mid(l_1,l_2)} \mu(l) \notag\\
  &= \sum_{dl \le x} d\mu(l) \Big( \sum_{mk \le x/(dl)} \mu(k) \Big)^2 \notag\\
  &= \sum_{dl \le x} d\mu(l) = \frac {x^2}{2\zeta(2)} + O(x\log x), \notag
\end{align}
on using variants of \eqref{e3} and \eqref{e4}. The asymptotic formula \eqref{eq1} follows from \eqref{eq3.1}--\eqref{eq3.3}. 

\section{The second moment: An analytic approach}
\label{s4}

In this section, we assume that $x$ and $y$ are as in part (ii) of Theorem~\ref{th2}. Recall the identity
\begin{equation}\label{eq4.1}
  \sum_{q = 1}^{\infty} \frac {c_q(n)}{q^s} = \frac {\sigma_{1-s}(n)}{\zeta(s)} \qquad (\re(s) > 1),
\end{equation}
where $\zeta(s)$ denotes the Riemann zeta-function and $\sigma_z(n) = \sum_{d \mid n} d^z$. We note that without loss of generality, we may assume that $x, y \in \mathbb Z + 1/2$. Then, using \eqref{eq1} and the truncated Perron formula \cite[Corollary 5.3]{MoVa07}, we get
\begin{equation}\label{eq2}
  \sum_{q \le x} c_q(n) = \frac 1{2\pi i} \int_{\alpha - iT}^{\alpha + iT} \frac {\sigma_{1-s}(n)}{\zeta(s)} \frac {x^s}s \, ds + E_1(x,n).
\end{equation}
Here, $T$ is a real parameter at our disposal, $\alpha > 1 + (\log x)^{-1}$ is a real number, and the remainder $E_1(x,n)$ satisfies
\[
  E_1(x,n) \ll \frac {x^{\alpha}}T \sum_{q = 1}^{\infty} \frac {|c_q(n)|}{q^{\alpha}} + \sum_{x/2 < q < 2x} |c_q(n)| \min\bigg( 1, \frac x{T|x - q|} \bigg).
\]
If we assume further that $\alpha < 1 + 3(\log x)^{-1}$, we can show that
\begin{equation}\label{eq3}
  E_1(x,n) \ll (x^2/T)\sigma_0(n)L,
\end{equation}
where $L = \log(Txy)$. Applying \eqref{eq2} and \eqref{eq3} with $\alpha = \alpha_j = 1 + j(\log x)^{-1}$, $j = 1, 2$, we obtain
\begin{equation}\label{eq4}
  \Big( \sum_{q \le x} c_q(n) \Big)^2 = \frac 1{(2\pi i)^2} \int_{\alpha_1 - iT}^{\alpha_1 + iT} \int_{\alpha_2 - iT}^{\alpha_2 + iT} F(s_1, s_2, n) \, ds_2ds_1 + E_2(x,n),
\end{equation}
where
\[
  F(s_1, s_2, n) = \frac {\sigma_{1-s_1}(n)\sigma_{1-s_2}(n)}{\zeta(s_1)\zeta(s_2)} \frac {x^{s_1+s_2}}{s_1s_2}
\]
and
\[
  E_2(x,n) \ll (x^4/T)\sigma_0(n)^2L^3 \big( x^{-1} + T^{-1} \big) \ll (x^3/T)\sigma_0(n)^2L^3,
\]
since we will later choose $T \ge x$. Summing \eqref{eq4} over $n$, we deduce
\begin{equation}\label{eq5}
  C_2(x,y) = \frac 1{(2\pi i)^2} \int_{\alpha_1 - iT}^{\alpha_1 + iT}
  \int_{\alpha_2 - iT}^{\alpha_2 + iT} \frac {G(s_1, s_2; y)}{\zeta(s_1)\zeta(s_2)}
  \frac {x^{s_1+s_2}}{s_1s_2} \, ds_2ds_1 + O\bigg( \frac {yx^3L^6}T \bigg),
\end{equation}
where 
\[
  G(s_1, s_2; y) = \sum_{n \le y} \sigma_{1-s_1}(n)\sigma_{1-s_2}(n).
\]

We now recall Ramanujan's identity
\begin{equation}\label{eq6}
  \sum_{n = 1}^{\infty} \frac {\sigma_a(n)\sigma_b(n)}{n^s} 
  = \frac {\zeta(s)\zeta(s-a)\zeta(s-b)\zeta(s-a-b)}{\zeta(2s-a-b)},
\end{equation}
valid for $\re(s) > \max\{ 1, 1 + \re(a), 1 + \re(b), 1 + \re(a+b) \}$. From \eqref{eq6}, by another application of Perron's formula, we get
\begin{equation}\label{eq7}
  G(s_1, s_2; y) = \frac 1{2\pi i} \int_{\alpha - 3iT}^{\alpha + 3iT} H(s_1, s_2, w) \frac {y^w}w \, dw + O\bigg( \frac {yL^4}T \bigg),
\end{equation}
where $\alpha = 1 + 3(\log y)^{-1}$, $y \in \mathbb Z + 1/2$, and
\[
  H(s_1,s_2,w) = \frac {\zeta(w)\zeta(w+s_1-1)\zeta(w+s_2-1)\zeta(w+s_1+s_2-2)}{\zeta(2w+s_1+s_2-2)}.
\]
Let $\Gamma(\alpha, \beta, T)$ denote the contour consisting of the line segments $[\alpha - iT, \beta - iT]$, $[\beta - iT, \beta + iT]$ and $[\beta + iT, \alpha + iT]$. We now move the integration in \eqref{eq7} to $\Gamma(\alpha, 1/2, 3T)$; let us denote the respective integrals $I_1, I_2$ and $I_3$. 

When $\re(s_j) = \alpha_j$, the integrals $I_1$ and $I_3$ over the two horizontal line segments are bounded above by
\begin{equation}\label{eq8}
  I_1, I_3 \ll \zeta(\alpha_1 + \alpha_2 - 1)L^4 \int_{1/2}^{\alpha} T^{1-2\sigma}y^{\sigma} d\sigma \ll L^5\big( y^{1/2} + yT^{-1} \big).
\end{equation}
Here, we have used the standard convexity bound (see \cite[eqn. (5.1.4)]{Titc86} for a slightly weaker version)
\begin{equation}\label{eq9}
  |\zeta(\sigma + it)| \ll (|t|+2)^{(1-\sigma)/2}\log(|t|+2) \qquad (1/2 \le \sigma \le 1).
\end{equation}
Furthermore, by H\"older's inequality,
\begin{equation}\label{eq10}
  I_2 \ll y^{1/2}\zeta(\alpha_1 + \alpha_2 - 1)\prod_{0 \le j,k \le 1} M_4(1/2 + j(\alpha_1 - 1) + k(\alpha_2-1), 5T),
\end{equation}
where for $1/2 \le \sigma \le 1$,
\[
  M_4^4(\sigma, T) = \int_{-T}^T |\zeta(\sigma + it)|^4 \, \frac {dt}{1 + |t|}.
\]
Appealing to the fourth-moment estimates for $\zeta(s)$ \cite[\S7.5 \& \S7.6]{Titc86}, we deduce from \eqref{eq10} that
\begin{equation}\label{eq11}
  I_2 \ll y^{1/2}L^6.
\end{equation}
Combining \eqref{eq7}, \eqref{eq8} and \eqref{eq11}, we obtain
\begin{equation}\label{eq12}
  G(s_1, s_2; y) = \sum_{j = 1}^4 R_j(s_1,s_2; y) + O( L^6( y^{1/2} + yT^{-1} ) ),
\end{equation}
where $R_j(s_1,s_2; y)$ are the four residues of the integrand in \eqref{eq7}:
\begin{align*}
  R_1(s_1,s_2; y) &= \frac {y\zeta(s_1)\zeta(s_2)\zeta(s_1+s_2-1)}{\zeta(s_1+s_2)}, \\
  R_2(s_1,s_2; y) &= \frac {y^{2-s_1}\zeta(2-s_1)\zeta(1 - s_1 + s_2)\zeta(s_2)}{(2-s_1)\zeta(2-s_1+s_2)}, \\
  R_3(s_1,s_2; y) &= \frac {y^{2-s_2}\zeta(2-s_2)\zeta(1 + s_1 - s_2)\zeta(s_1)}{(2-s_2)\zeta(2+s_1-s_2)}, \\
  R_4(s_1,s_2; y) &= \frac {y^{3-s_1-s_2}\zeta(3-s_1-s_2)\zeta(2-s_2)\zeta(2-s_1)}{(3-s_1-s_2)\zeta(4-s_1-s_2)}.
\end{align*}
Substituting \eqref{eq12} into the right side of \eqref{eq5}, we get
\begin{equation}\label{eq13}
  C_2(x,y) = \sum_{j = 1}^4 C_{2,j}(x,y) + O( yx^2L^{10}( y^{-1/2} + xT^{-1}) ),
\end{equation}
where
\[
  C_{2,j}(x,y) = \frac 1{(2\pi i)^2} \int_{\alpha_1 - iT}^{\alpha_1 + iT} \int_{\alpha_2 - iT}^{\alpha_2 + iT} \frac {R_j(s_1, s_2; y)}{\zeta(s_1)\zeta(s_2)} \frac {x^{s_1+s_2}}{s_1s_2} \, ds_2ds_1.
\]

We proceed to evaluate the integrals $C_{2,j}(x,y)$, starting with $C_{2,1}(x,y)$. We first move the integration over $s_2$ to the contour $\Gamma(\alpha_2, 1/2, T)$; let us denote the integrals over the three line segments $J_{1,1}, J_{1,2}$ and $J_{1,3}$. The contributions from the integrals over the two horizontal line segments are bounded using \eqref{eq9}:
\[
  J_{1,1}, J_{1,3} \ll \frac {xyL}T  \int_{-T}^T \frac {dt}{1 + |t|} \int_{1/2}^{\alpha_2} T^{(1-\sigma)/2}x^{\sigma} d\sigma \ll xyL^2( xT^{-1} + 1 ).
\]
Furthermore, using estimates for the mean square of $\zeta(s)$ in the critical strip, we find that the contribution from the integral over the line $\re(s_2) = 1/2$ is
\begin{align*}
  J_{1,2} &\ll yx^{3/2} \int_{-T}^{T} \int_{-T}^{T} \frac {|\zeta(\alpha_1-1/2 + i(t_1+t_2))|} {(|t_1|+1)(|t_2| + 1)} \, dt_1dt_2 \\
  &\ll yx^{3/2} \int_{-2T}^{2T} |\zeta(\alpha_1-1/2 + iu)| \int_{-T}^T \frac {dtdu}{(|t|+1)(|t-u| + 1)} \\
  &\ll yx^{3/2}L \int_{-2T}^{2T} |\zeta(\alpha_1-1/2 + iu)| \frac {du}{|u|+1} \ll yx^{3/2}L^3.
\end{align*}
Hence, accounting for the residue at $s_2 = 2 - s_1$, we have
\begin{align}\label{eq14}
  C_{2,1}(x,y) &= \frac {yx^2}{\zeta(2)}\frac 1{2\pi i} \int_{\alpha_1 - iT}^{\alpha_1 + iT} \frac {ds_1}{s_1(2-s_1)} + O( yx^2L^3( x^{-1/2} + T^{-1}) ), \\
  &= \frac {yx^2}{2\zeta(2)} + O( yx^2L^3( x^{-1/2} + T^{-1}) ). \notag
\end{align}

To estimate $C_{2,4}(x,y)$, we move the integration over $s_2$ to the contour $\Gamma(\alpha_2, \beta, T)$, where $\beta = 5/2 - \alpha_1$; let us denote the respective integrals $J_{4,1}, J_{4,2}$ and $J_{4,3}$. Similarly to the estimation of $J_{1,j}$, we have
\[
  J_{4,1}, J_{4,3} \ll yx^2L^5 ( T^{-1} + (x/yT)^{1/2} ), \quad J_{4,2} \ll y^{1/2}x^{5/2}L^5.
\]
Since the integrand is holomorphic in the region between the two contours, we deduce that
\begin{align}\label{eq15}
  C_{2,4}(x,y) &\ll yx^2L^5 ( (x/y)^{1/2} + T^{-1} ).
\end{align}
Similarly, by moving the integration to $\Gamma(\alpha_2, 3/2, T)$, we find that
\begin{align}\label{eq16}
  C_{2,3}(x,y) &\ll yx^2L^5 ( (x/y)^{1/2} + T^{-1} ),
\end{align}
and
\begin{align}\label{eq17}
  C_{2,2}(x,y) &= \frac {yx^2}{\zeta(2)} K(x,y) + O( yx^2L^5 ( (x/y)^{1/2} + T^{-1} ) ),
\end{align}
where
\[
  K(x,y) = \frac 1{2\pi i} \int_{\alpha_1 - iT}^{\alpha_1 + iT} \frac {\zeta(2-s)}{\zeta(s)} \frac {(yx^{-2})^{1-s}}{s^2(2-s)} \, ds.
\]
The first term on the right side of \eqref{eq17} arises from the pole of $R_2(s_1, s_2; y)$ at $s_2 = s_1$. Thus, by \eqref{eq13}--\eqref{eq17},
\begin{equation}\label{eq18}
  C_2(x,y) = \frac {yx^2}{\zeta(2)}\big( 1/2 + K(x,y) \big) + O( yx^2L^{10}( x^{-1/2} + (x/y)^{1/2} + xT^{-1}) ).
\end{equation}
It remains to evaluate the integral $K(x,y)$. 

Recall that $yx^{-2} \le L^B$. We move the integration to the contour $\Gamma(\alpha_1, 1, T)$ and denote the integrals along the edges $K_1(x,y)$, $K_2(x,y)$ and $K_3(x,y)$. Then, by the bounds for $\zeta(s)$ in \cite[\S6.19]{Titc86}, we have
\[
  K_1(x,y), K_3(x,y) \ll T^{-2} \quad \text{and} \quad
  K_2(x,y) = \kappa(u) + O( T^{-2} ),
\]
where $u = \log(yx^{-2})$ and $\kappa(u)$ is the Fourier integral
\begin{equation}\label{eq22}
  \kappa(u) = \frac 1{2\pi} \int_{-\infty}^\infty f(it) e^{-itu} \, dt, \qquad
  f(s) = \dfrac{\zeta(1-s)}{\zeta(1+s)} \frac{1}{(1 + s)^2 (1 - s)}.
\end{equation}
We choose $T = x^2$ and get the asymptotic formula
\begin{equation}\label{eq20}
  C_2(x,y) = \frac {yx^2}{\zeta(2)}( 1/2 + \kappa(u) ) + O( yx^2L^{10}( (x/y)^{1/2} + x^{-1/2} ) ).
\end{equation}
Thus, to complete the proof of part (ii) of the theorem, we need only establish the desired properties of $\kappa(u)$. We examine $\kappa(u)$ in the next section.

\section{Estimation of $\kappa(u)$}
\label{s5}

Suppose first that $|u| \to \infty$. For a fixed $\epsilon > 0$ and a real $t$, we define
\[
  \delta(t) = C_{\epsilon}(\log(|t| + 2))^{-2/3-\epsilon}.
\]
By the standard bounds for $\zeta(s)$ near the edge of the critical strip \cite[p.~135]{Titc86},
\[
  |f(\sigma + it)| \ll \frac {\delta(t)^{-2}}{(1 + |t|)^3} \quad \text{when } |\sigma| \le \delta(t).
\]
Applying Cauchy's integral formula for the circle $\gamma_t$, $|s - it| = \delta(t)$, we deduce that
\[
  |f^{(k)}(it)| = \bigg| \frac{k!}{2 \pi i} \int_{\gamma_t} \frac{f(s)}{(s - it)^{k+1}} \, ds \bigg| \ll \frac{k!\delta(t)^{-k-2}}{(1 + |t|)^3} \quad (k \ge 1).
\]
Hence, a $k$-fold integration by parts gives
\begin{align*}
  \kappa(u) &\ll \frac {1}{|u|^{k}} \bigg| \int_{-\infty}^\infty f^{(k)}(it)e^{-itu} \, dt \bigg|
  \ll \frac {k!C_{\epsilon}^{-k-2}}{|u|^{k}} \int_2^{\infty} (\log t)^{(2/3 + \epsilon)(k + 2)}t^{-3} \, dt.
\end{align*}
Using the inequality
\[
  \int_2^{\infty} (\log t)^{\alpha}t^{-3} \, dt \le \int_0^{\infty} z^{\alpha}e^{-2z} \, dz = \frac {\Gamma(\alpha + 1)}{2^{\alpha+1}},
\]
we deduce that
\begin{align*}
  \kappa(u) &\ll \frac{k!m_k!}{2^{m_k}(C_{\epsilon}|u|)^k} \ll \frac{k^{(5/3 + \epsilon)k}}{(C_{\epsilon}|u|)^k},
\end{align*}
where $m_k = \lceil (2/3 + \epsilon)(k + 2) \rceil$. Finally, we set $k = \lfloor C_{\epsilon}'|u|^{1/(5/3+\epsilon)} \rfloor$, and we get
\[
  \kappa(u) \ll \exp( -|u|^{3/5 - \epsilon} ).
\]

Next, we focus on the case when $|u|$ is bounded. The integral $\kappa(u)$ is differentiable, and
\[
  |\kappa'(u)| \le \frac 1{2\pi} \int_{-\infty}^{\infty} \frac {|t| \, dt}{(1 + t^2)^{3/2}} \le \frac 12.
\]
Hence,
\begin{equation}\label{eq25}
  |\kappa(u) - \kappa(v)| \le \frac {|u - v|}2.
\end{equation}
Since the integrand is bounded above by $|t|^{-3}$, we have
\begin{equation}\label{eq26}
  \kappa(u) = \frac 1{2\pi} \int_{-60}^{60} f(it)e^{-itu} \, dt + \vartheta_u 10^{-4},
\end{equation}
where $|\vartheta_u| \le 1/2$. Let $\kappa_0(u)$ denote the integral on the right side of \eqref{eq26}. By partial integration,
\[
  |\kappa_0(u)| \le \frac 1{2\pi|u|} \Big( |f(60i)| + |f(-60i)| + \int_{-60}^{60} |f'(it)| \, dt \Big) \le \frac {4.2}{2\pi|u|} \le 0.67|u|^{-1}.
\]
Combining the latter inequality and \eqref{eq26}, we deduce that $|\kappa(u)| \le 0.395$ when $|u| \ge 1.7$. It remains to estimate $\kappa(u)$ when $|u| \le 1.7$. 

In view of \eqref{eq25} and \eqref{eq26}, it suffices to tabulate $\kappa_0(u)$ for a set of suitably chosen values $u_1, u_2, \dots$, with $|u_j| \le 1.7$. To this end, a quick calculation using \emph{Mathematica} yields the following table:
\begin{center}
\medskip
\begin{tabular}{ccl} \midrule
$u_j$ & $\kappa_0(u_j)$ & \hfill Covered range \hfill \\ \midrule
$-1$  & $-0.0277$ & $-1.7 \le u \le -0.3$ \\
$0.3$  & $-0.0984$ & $-0.3 \le u \le \phantom{+}0.9$ \\
$1.1$ & $-0.2586$ & $\phantom{+}0.9 \le u \le \phantom{+}1.3$ \\
$1.5$ & $-0.2925$ & $\phantom{+}1.3 \le u \le \phantom{+}1.7$ \\ \midrule
\end{tabular}
\medskip
\end{center}
Each row of this table lists a value of $\kappa_0(u_j)$ and a range of the variable $u$, where we can infer the bound $\kappa(u) \ge -0.4$ from \eqref{eq25}, \eqref{eq26} and the listed value of $\kappa_0(u_j)$. For example, when $0.9 \le u \le 1.3$, we have
\[
  \kappa(u) \ge \kappa(1.1) - 0.5|u - 1.1| > \kappa_0(1.1) - 0.1001 > -0.36.
\]
Clearly, this completes the proof of part (ii) of Theorem \ref{th2}.\\

\paragraph*{\em Acknowledgment} The authors would like to thank the referee for his/her comments, which helped improve the presentation.


\begin{thebibliography}{9}
  \bibitem {ChKu05} {T. H. Chan and A. V. Kumchev},
  \emph{Approximating reals by sums of rationals}, 
  unpublished.
  \bibitem {GrKo91} {S. W. Graham and G. Kolesnik},
  \emph{Van der Corput's Method of Exponential Sums},
  Cambridge Univ. Press, 1991.
  \bibitem {MoVa07} {H. L. Montgomery and R. C. Vaughan},
  \emph{Multiplicative Number Theory I. Classical Theory},
  Cambridge Univ. Press, 2007.
  \bibitem {Titc86} {E. C. Titchmarsh},
  \emph{The Theory of the Riemann Zeta-Function},
  2nd ed., Oxford Univ. Press, 1986.
\end{thebibliography}
\end{document}